\documentclass[11pt]{amsart}

\usepackage{amsmath}

\usepackage{amssymb}

\usepackage{amscd}

\usepackage{xspace}

\usepackage{verbatim}

\newcommand{\fnz}{\footnotesize}

\newcommand{\be}{\begin{equation}}

\newcommand{\bel}[1]{\begin{equation}\label{#1}}

\newcommand{\ee}{\end{equation}}%% This macro does not work with amstex.

\newtheorem{subn}{\name}

\newcommand{\bsn}[1]{\def\name{#1}\begin{subn}}

\newcommand{\esn}{\end{subn}}

\newtheorem{sub}{\name}%[section]

\newcommand{\bs}{\begin{sub}}

\newcommand{\es}{\end{sub}}

\newcommand{\bsl}[1]{\begin{sub}\label{#1}}

\newcommand{\bth}[1]{\def\name{Theorem}\begin{sub}\label{t:#1}}

\newcommand{\blemma}[1]{\def\name{Lemma}\begin{sub}\label{l:#1}}

\newcommand{\bcor}[1]{\def\name{Corollary}\begin{sub}\label{c:#1}}

\newcommand{\bdef}[1]{\def\name{Definition}\begin{sub}\label{d:#1}}

\newcommand{\bprop}[1]{\def\name{Proposition}\begin{sub}\label{p:#1}}

\newcommand{\BA}{\begin{array}}

\newcommand{\EA}{\end{array}}

\newcommand{\BAN}{\renewcommand{\arraystretch}{1.2}

\setlength{\arraycolsep}{2pt}\begin{array}}

\newcommand{\BAV}[2]{\renewcommand{\arraystretch}{#1}

\setlength{\arraycolsep}{#2}\begin{array}}

\newcommand{\BSA}{\begin{subarray}}

\newcommand{\ESA}{\end{subarray}}

\newcommand{\BAL}{\begin{aligned}}

\newcommand{\EAL}{\end{aligned}}

\newcommand{\BALG}{\begin{alignat}}

\newcommand{\EALG}{\end{alignat}}%% the abbrev. does not work with latex2e

\newcommand{\BALGN}{\begin{alignat*}}

\newcommand{\EALGN}{\end{alignat*}}%% the abbrev. does not work with latex2e

\newcommand{\note}[1]{\textit{#1.}\hspace{2mm}}

\newcommand{\Remark}{\note{Remark}}

\newcommand{\forevery}{\quad \forall}

\newcommand{\abs}[1]{\left |#1\right |}%% adjustable vertical delimiters

\newcommand{\myfrac}[2]{{\displaystyle \frac{#1}{#2} }}

\newcommand{\myint}[2]{{\displaystyle \int_{#1}^{#2}}}

\newcommand{\prt}{\partial}

\newcommand{\ti}{\times}

\newcommand{\nind}{\noindent}

            \def\gg{\gamma}

       \def\gd{\delta}      \def\ge{\epsilon}

            \def\gl{\lambda}

\def\gm{\mu}

\def\gs{\sigma}       \def\gt{\tau}

\def\Gg{\Gamma}     \def\Gd{\Delta}

    \def\Gs{\Sigma}

   \def\CB{{\mathcal B}}

   \def\BBR {\mathbb R}

\begin{document}%

\nind{\small {\bf \'Equations aux d\'eriv\'ees partielles}}/{\small {\it Partial Differential Equations}}\\

\nind{\LARGE\bf Sym\'etrie des grandes solutions d'\'equations elliptiques semi lin\'eaires 
\footnote {A para\^\i tre dans{ \it C. R. Acad. Sci. Paris, Ser. I.}}}\\

\nind{\bf Alessio Porretta}
\footnote{L'auteur a b\'en\'efici\'e du support du projet europ\'een RTN:
FRONTS-SINGULARITIES, RTN contract: HPRN-CT-2002-00274. }\\
{\small Dipartimento di Matematica, 
Universit\`a di Roma "Tor Vergata", Via della Ricerca Scientifica 1, 00133 Roma, Italia\\
 Email: {\it porretta@mat.uniroma2.it }} \\[-2mm]

\nind {\bf Laurent V\'eron}\\
{\small Laboratoire de Math\'ematiques et Physique Th\'eorique, CNRS UMR 6083, Facult\'e des Sciences, 37200 Tours, France. \\
Email: {\it veronl@univ-tours.fr}}\\
\underline{\phantom{-----------------------------------------------------------------------------------------------}}

\begin{quotation}
{\fnz{\sc \bf{R\'esum\'e}.}
Soit $g$ une fonction localement lipschitzienne de la variable r\'eelle. On suppose que $g$ v\'erifie la condition de Keller et Osserman et qu'il existe un r\'eel $a>0$  tel que $g$ est convexe sur  $[a,+\infty[$. Alors toute solution $u$ de $-\Gd u+g(u)=0$ dans une boule $B$ de $\BBR^N$, $N\geq 2$, qui tend vers l'infini au bord de $B$, est une fonction radiale. 
 }\end {quotation}

\indent {\large\bf Symmetry of large solutions of semilinear elliptic equations}
\begin{quotation}
 {\fnz{\sc \bf{Abstract}.}
 \setcounter{equation}{0}
Let $g$ be a locally Lipschitz continuous function defined on $\BBR$. We assume that $g$ satisfies the Keller-Osserman condition and there exists a positive real number $a$ such that $g$ is convex on $[a,\infty)$. Then any solution $u$ of $-\Gd u+g(u)=0$ in a ball $B$ of $\BBR^N$, $N\geq 2$, which tends to infinity on $\prt B$, is spherically symmetric.
 }
\end {quotation}\underline{\phantom{-----------------------------------------------------------------------------------------------}}\\

\nind{\large\bf{Abridged English version}}\\

\nind Let $g: \BBR\mapsto\BBR$ be a locally Lipschitz continuous function and $B_{R}(0)$ the open $N$-ball ($N\geq 2$) of center $0$ and radius $R>0$. A classical result due to Gidas, Ni and Nirenberg \cite {GNN} asserts that any positive solution $u$ of
\begin {equation}\label {equ}
-\Gd u+g(u)=0
\end {equation}
in $B_{R}(0)$ which vanishes on $\prt B_{R}(0)$ is radial. A conjecture proposed by H. Brezis is that any {\it large solution} of (\ref {equ}), that is a solution which verifies
\begin {equation}\label {b-u}
\lim_{\abs x\to R} u(x)=\infty, 
\end {equation}
is radial. The existence of such solution is ensured by the Keller-Osserman condition: there exists some $a>0$ such that $g$ is nondecreasing on $[a,\infty)$ and
\begin {equation}\label {KO}
\myint{a}{\infty}\myfrac {ds}{\sqrt {G(s)}}<\infty\quad\mbox {where }\;G(s)=\myint{a}{s} g(\gs)d\gs.
\end {equation}
We prove two symmetry  results dealing with this conjecture.

 %%%%%%%%%%THEOREM 1%%%%%%%%%%%%%%%%%%%%%%
\bth{main1} Assume $g$ is locally Lipschitz continuous and let $u$ be a large solution of (\ref {equ}) in a ball $B=B_R(0)\subset \BBR^N$, $N\geq 2$. If there holds
\begin{equation}\BA{l}\label {cond}
(i)\quad\lim_{\abs x\to R}\myfrac{\prt u}{\prt r}(x)=\infty\\
(ii) \quad\abs{\nabla_\gt u(x)}=\circ\left(\myfrac{\prt u}{\prt r}(x)\right)\quad \mbox {as }
\abs x\to R,
\EA\end{equation}
then $u$ is radial and $\myfrac{\prt u}{\prt r}(x)>0$ on $B_R(0)\setminus\{0\}$.
\es

%%%END THEOREM 1%%%%%%%%%%%%%%%%%%%%%%%%%%%%%%
In this statement 
$\myfrac{\prt u}{\prt r}(x)=\langle Du(x),x/\abs x\rangle$ is the radial derivative and 
$\nabla_\gt u(x)= Du(x)-{\abs x}^{-2}\langle Du(x),x\rangle\, x$ is the tangential gradient. This result is settled upon an adaption of the key lemma of \cite{GNN} in the framework of large solutions. Next we give a sufficient condition in order (\ref{cond}) to hold.
 %%%%%%%%%%THEOREM 2%%%%%%%%%%%%%%%%%%%%%%
\bth{main2} Assume $g$ is  locally Lipschitz continuous, convex on $[a,+\infty)$ for some $a>0$ and satisfies (\ref{KO}). Then any large solution of (\ref {equ}) in a ball is a radial function.
\es
%%%END THEOREM 2%%%%%%%%%%%%%%%%%%%%%%%%%%%%%%

\nind\Remark It is important to notice that this result is not related with uniqueness. For example, if $g(x)=x^2$ it is known that uniqueness may not hold if the radius of the ball is large enough. As a striking example, if $g$ is any polynomial of degree larger than one with positive coefficient of higher order, any large solution of (\ref {equ}) in a ball is radial.

%%%%%%%%%%%%%%%%%%%%%%%%%%%%%%%%%%%%%%%%%%%%
%%%%%VERSION FRANCAISE%%%%%%%%%%%%%%%%%%%%%%%%%%%
\medskip

%%%%%%%%%%%%%%%%%%%%%%%%%%%%%%%%%%%%%%%%%%%%
\setcounter {equation}{0}
%%%%%%%%%%%%%%%%%%%%%%%%%%%%%%%%%%%%%%%%%%
%%%%%%%%%%%%%%%%%%%%%%%%%%%%%%%%%%%%%%%%%%%%
\nind{\large\bf{R\'esultats principaux}}\\ [-2mm]

Soit $g:\BBR\mapsto\BBR$ une fonction localement lipschitzienne et $B_R(0)$ la boule de centre $0$ et de rayon $R$ dans $\BBR^N$, $N\geq 2$. Un r\'esultat classique du \`a Gidas, Ni et Nirenberg \cite {GNN} affirme que si $u$ est une solution positive de 
\begin {equation}\label {equ}
-\Gd u+g(u)=0
\end {equation}
dans $B_R(0)$ qui s'annulle sur $\prt B_R(0)$ alors elle est radiale. Si $u$ prend la valeur $k$ au bord le r\'esultat reste valable pourvu que $u-k$ ne change pas de signe dans $B_R(0)$. Partant de cette observation, H. Brezis a conjectur\'e que si $u$ est une {\it grande solution}, c'est \`a dire une solution qui v\'erifie
\begin {equation}\label {ex}
\lim_{\abs x\to R}u(x)=\infty,
\end {equation}
alors elle est radiale.  L'existence de grandes solutions est associ\'ee \`a la condition de Keller et Osserman qui est satisfaite si $g$ est positive et croissante sur $[a,+\infty[$ pour un $a>0$ et y v\'erifie
\begin {equation}\label {KO}
\myint{a}{+\infty}\myfrac {ds}{\sqrt {G(s)}}<+\infty\quad\mbox {o\`u }\;G(s)=\myint{a}{s} g(\gs)d\gs.
\end {equation}
 Nous donnons deux r\'esultats qui confirment la validit\'e de la conjecture de Brezis.\\

 %%%%THEOREME 1%%%%%%%%%%%%%%%%%%%%%%%%%%%%%%

\nind{\bf Th\'eor\`eme 1} {\it Supposons que $g$ est localement lipschitzienne et soit $u$ une grande solution de (\ref {equ}) dans la boule $B=B_R(0)\subset \BBR^N$, $N\geq 2$. Si on a
\begin{equation}\BA{l}\label {cond}
(i)\quad\lim_{\abs x\to R}\myfrac{\prt u}{\prt r}(x)=\infty\\
(ii) \quad\abs{\nabla_\gt u(x)}=\circ\left(\myfrac{\prt u}{\prt r}(x)\right)\quad \mbox {as }
\abs x\to R,
\EA\end{equation}
alors $u$ est radiale et $\myfrac{\prt u}{\prt r}(x)>0$ dans $B_R(0)\setminus\{0\}$.}\\
%%%FIN THEOREME 1 %%%%%%%%%%%%%%%%%%%%%%%%%%%%%%%%

Dans cet \'enonc\'e $\myfrac{\prt u}{\prt r}(x)=\langle Du(x),x/\abs x\rangle$ est la d\'eriv\'ee radiale de $u$ et 
$\nabla_\gt u(x)= Du(x)-{\abs x}^{-2}\langle Du(x),x\rangle\, x$ son gradient tangentiel. \\ 

 %%%%THEOREME 2%%%%%%%%%%%%%%%%%%%%%%%%%%%%%%
\nind {\bf Th\'eor\`eme 2} {\it Supposons que $g$ est localement lipschitzienne et qu'il existe $a>0$ tel que $g$ est convexe sur $[a,+\infty[$ et y v\'erifie (\ref {KO}). Alors toute grande solution de (\ref {equ}) dans une boule est radiale.}\\
%%%FIN THEOREME 2 %%%%%%%%%%%%%%%%%%%%%%%%%%%%%%%%

\nind{\it Remarque.} Il est important de noter que ce r\'esultat n'augure en rien de l'unicit\'e des grandes solutions de (\ref {equ}). Ainsi, si $g(x)=x^2$, il est classique \cite {Po} que si le rayon de la boule est assez grand, il existe plusieurs grandes solutions, dont une seule positive. Par exemple, si $g$ est un polyn\^ome de degr\'e $>1$ dont le coefficient du terme de plus haut degr\'e est positif, alors le r\'esultat du Th\'eor\`eme 2 s'applique.\\

Le r\'esultat suivant \'etend aussi un autre th\'eor\`eme de \cite {GNN}.\\

\nind {\bf Corollaire 1} {\it Supposons que $g$ v\'erifie les hypoth\`eses du Th\'eor\`eme 2. Si $u$ est une solution de (\ref {equ}) dans $\Gg_{R,r}=\{x\in \BBR^N:r< \abs x<R\}$ qui v\'erifie (\ref {ex}), alors $\myfrac{\prt u}{\prt r}(x)>0$ pour tout $x\in \Gg_{R,(r+R)/2}$.
}\\

\nind {\it Principe de la d\'emonstration du Th\'eor\`eme 1}. On commence par  noter que pour tout $P\in\prt B^+=\prt B_{R}(0)\cap \{x_{1}>0\}$, il existe $\gd\in ]0,R[$ tel que 
\begin {equation}\label {est1}
\myfrac{\prt u}{\prt x_{1}}(x)>0\forevery x\in B_{R}(0)\cap B_{\gd}(P).
\end {equation}
Ceci d\'ecoule imm\'ediatement de (\ref{cond}). La suite de la d\'emonstration du Th\'eor\`eme 1 repose sur  la m\'ethode des plans mobiles comme dans \cite {GNN}. Soit
$\CB=\{{\bf e}_1,...{\bf e}_N\}$ une base orthonorm\'ee de $\BBR^N$ et $(x_1,...,x_N)$ les coordonn\'ees d'un point $x$ dans cette base. Pour $0<\gl<R$ on d\'esigne par $T_\gl$ l'hyperplan $\{x:x_1=\gl\}$, $\Gs_\gl=\{x\in B_R(0):\gl<x_1<R\}$, $\Gs'_\gl=\{x\in B_R(0):2\gl-R<x_1<\gl\}$,  par
$x_\gl$ le sym\'etrique de $x$, par rapport \`a $T_\gl$, de coordonn\'ees
$(2\gl-x_1,x_2,...,x_N)$ et par $u_{\gl}$ la fonction r\'efl\'echie de $u$, d\'efinie par $u_{\gl}(x)=u(x_{\gl})$. On applique (\ref{est1}) avec $P=P_{0}=R{\bf e}_1$, $\gd_{0}=\gd(P_{0})$. On en d\'eduit que pour tout $\gl\in [\gl_{0},R[$ (o\`u $\gl_{0}=R-\gd_{0}^2/2R$) on a
\begin {equation}\label {est2}
u(x_{\gl})<u(x)\quad\mbox {et }\;\myfrac{\prt u}{\prt x_{1}}(x)>0\forevery x\in \Gs_{\gl}.
\end {equation}
Soit $\gm=\inf \{\gl>0:\;\mbox {t. q. (\ref{est2}) soit v\'erifi\'ee}\}$. On suppose $\gm>0$. Par d\'efinition $u\geq u_{\gm}$ dans $\Gs_{\gm}.$ Soit $K_{\gm}=T_{\gm}\cap \prt B_{R}(0)$. Comme $K_{\gm}$ est compact, grace \`a (\ref {est1}) il existe un $\ge$-voisinage $U_{\ge}$ de $K_{\gm}$ tel que 
\begin {equation}\label {est3}\myfrac{\prt u}{\prt x_{1}}(x)>0\forevery x\in U_{\ge}\cap B_{R}(0).
\end {equation}
On pose $D_{\ge}=B_{R-\ge/2}(0)\cap \Gs_{\gm}$ et $a(x)=(g(u)-g(u_\mu)/(u-u_\mu)$.  Comme $w=u-u_{\gm}$ v\'erifie
\begin {equation}\label {est4}
\Gd w-aw=0\quad \mbox {dans }D_{\ge},\quad w\geq 0, \quad w\equiv \!\!\!\!\!/ \;\,0,
\end {equation}
on en d\'eduit $w>0$ par le principe du maximum fort, et $\prt u/\prt x_{1}>0$ sur $T_{\gm}\cap \prt D_{\ge}$ par le lemme de Hopf. La continuit\'e de $Du$  \`a l'int\'erieur et (\ref {est3}) impliquent qu'il existe $\gs>0$ tel que
\begin {equation}\label {est5}
\myfrac{\prt u}{\prt x_{1}}(x)>0\forevery x\in B_{R}(0)\cap\{x:\gm-\gs<x_{1}<\gm+\gs\}.
\end {equation}
De plus, comme $\ge$ est arbitrairement petit, $u>u_{\gm}$ dans $\Sigma_\mu$. La d\'efinition de $\gm$ implique qu'il existe une suite positive croissante $\{\gl_{n}\}$ convergeant vers $\gm$ et une suite de points 
$\{x_{n}\}$ convergeant vers  $\bar x\in \overline {\Gs_{\gm}}$ telles que 
$u(x_{n})\leq u((x_{n})_{\gl_{n}})$. Comme $u>u_{\gm}$ dans $\Sigma_\mu$, $\bar x$ ne peut appartenir \`a 
$\Gs_{\gm}$. Le th\'eor\`eme des accroissement fini et (\ref {est5}) impliquent que $\bar x$ ne peut appartenir non plus \`a $T_\gm$. Enfin $\bar x$ ne peut appartenir \`a $\Gs_{\gm}\setminus T_{\gm}$ puisque cela impliquerait que  $u(x_{n})-u((x_{n})_{\gl_{n}})$ tende vers $+\infty$. Par contradiction il s'ensuit que $\gm=0$. Changeant $x_{1}$ en $-x_{1}$ puis permutant les directions, on en d\'eduit que $u$ est radiale.\smallskip

 \nind {\it Principe de la d\'emonstration du Th\'eor\`eme 2}. La clef est le r\'esultat suivant.\\

\nind{\bf Lemme 1} {\it 
Supposons que $g$ v\'erifie les hypoth\`eses du Th\'eor\`eme 2, et que $u$ est une grande solution de (\ref {equ}) dans $B_{R}(0)$. Alors
\begin{equation}\label{infde}
\BA{l} 
(i)\quad\displaystyle\lim_{|x|\to R}\nabla_{\tau}u(x)= 0\\ 
[3mm]

(ii)\quad\displaystyle\lim_{|x|\to R}\myfrac{\partial u}{\partial r}(x)=\infty,
\EA
\end{equation}
et les deux limites ont lieu uniform\'ement par rapport \`a $\{x:\abs x=r\}$.}\\

\nind {\it D\'emonstration}. Soient $(r,\gs)\in\BBR_{+}\ti S^{N-1}$ les coordonn\'ees sph\'eriques dans $\BBR^N$, $\tilde\gs\in S^{N-1}$ et $\{\gg_{j}\}_{j=1}^{N-1}$ un ensemble de g\'eod\'esiques de
$S^{N-1}$ se coupant orthogonalement en $\tilde\gs$, par exemple $\gg_{j}(t)=e^{tA_{j}}(\tilde\gs)$ o\`u les matrices $\{A_{j}\}_{j=1}^{N-1}$ sont anti-sym\'etriques et v\'erifient $\langle A_{j}\tilde\gs,A_{k}\tilde\gs\rangle=\gd_{j}^k$. Si $\Gd_{S}$ est l'op\'erateur de Laplace-Beltrami sur $S^{N-1}$, on a
\begin {equation}\label {LB}
\Gd_{S}u(r,\tilde\gs)=\sum_{j=1}^{N-1}\myfrac{d^2u(r,\gg_{j}(t))}{dt^2}|_{t=0}.
\end {equation}
Par hypoth\`ese $g=g_{\infty}+\tilde g$ o\`u $g_{\infty}$ est convexe et v\'erifie (\ref {KO}) et 
$\tilde g$ est localement lipschitzien et identiquement nul sur $[M,+\infty[$ pour un $M>0$. Sans restriction on peut supposer $g_{\infty}$ croissante. Il existe $r_{0}\in ]0,R[$ tel que 
$u(x)\geq M$ pour tout $\abs x\geq r_{0}$. Ainsi
\begin {equation}\label {est51}
\abs {\Gd u-g_{\infty}(u)}=|\tilde g(u)|=|\tilde g(u)\chi_{_{B_{r_{\tiny{0}}}(0)}}|\leq K_{0}.
\end {equation}
Soit $\phi(x)=(2N)^{-1}(R^2-\abs x^2)$. Comme $\Gd\phi=-1$ on d\'eduit de (\ref{est51})
$$\Gd(u-K_{0}\phi)\geq g_{\infty}(u)\geq g_{\infty}(u-K_{0}\phi),
$$
et donc $u-K_{0}\phi$ est une sous-solution du probl\`eme
\begin {equation}\label {est6}\left\{\BA {l}
-\Gd v+g_{\infty}(v)=0\quad\mbox {in }\;B_{R}(0)\\
\;\lim_{\abs x\to R}v(x)=\infty.
\EA\right.\end {equation}
Par convexit\'e (voir par exemple \cite {MV}, \cite {MV1}, m\^eme si il existe une 
d\'emonstration plus directe dans le cas radial) ce probl\`eme admet une unique solution $v=U_{R}$. Comme $u+K_{0}\phi$ est une sur-solution, on en d\'eduit
\begin {equation}\label {est7}
U_{R}-K_{0}\phi\leq u\leq U_{R}+K_{0}\phi.
\end {equation}
Soit $h>0$, $j=1,...,N-1$ et $u^{h}(x)=u(e^{hA_{j}}(x))=u(r,e^{hA_{j}}\gs)$, o\`u $x=(r,\gs)$. Comme le probl\`eme est invariant par rotation, $u^h$ v\'erifie aussi (\ref{est7}). Par suite 
\begin {equation}\label {est8}
\lim_{\abs x\to R}u(x)-u^h(x)=0.
\end {equation}
De plus $\Gd u^h=g_{\infty}(u^h)$ dans 
$\Gg_{R,r_{0}}=B_{R}(0)\setminus B_{r_{0}}(0)$ et il existe $L>0$, ind\'ependant de $h$, tel que $\abs {(u-u^h)(x)}\leq L\abs h$ pour $\abs x=r_{0}$. Si $\Psi$ est la fonction harmonique dans $\Gg_{R,r_{0}}$, nulle sur 
$\prt B_{R}(0)$ et valant $1$ sur $\prt B_{r_{0}}(0)$, et $v^h=u^h+\abs hL\Psi$, alors
\begin {equation}\label {est9}
\Gd(v^h-u)\leq g_{\infty}(v^h)-g_{\infty}(u)\quad \mbox {dans }\;\Gg_{R,r_{0}}.
\end {equation}
La relation (\ref {est8}) implique que $v^h(x)-u(x)\to 0$ si $\abs x\to R$. Par monotonie 
$v^h=u^h+\abs hL\Psi\geq u$. Si on d\'efinit la d\'eriv\'ee de Lie selon le champ de vecteurs $\gs\mapsto A_{j}\gs$ par 
$$L_{A_{j}}u(r,\gs)=\myfrac {du(r,e^{tA_{j}\gs})}{dt}|_{t=0},$$
alors
\begin {equation}\label {est10}
|L_{A_{j}}u(r,\tilde\gs)|\leq L\Psi(x)\leq C(R-r).
\end {equation}
Cette relation implique (\ref{infde})-i. \smallskip

\nind Pour d\'emontrer (\ref{infde})-ii, on pose $w^h=h^{-2}(u^h+u^{-h}-2u)$. La convexit\'e de $g_{\infty}$ implique que $w^h$ 
v\'erifie $\Gd w^h\geq \xi(x)w^h$ dans $\Gg_{R,r_{0}}$,
o\`u $\xi(x)\geq 0$,  et donc que $w^h_{+}$ est sous-harmonique dans  $\Gg_{R,r_{0}}$. Comme $u$ est de classe $C^2$, il existe $\tilde L>0$, ind\'ependant de $h$, tel que $w^h\leq \tilde L$ sur $\prt B_{r_{0}}(0)$. Comme $w^h$ et $\Psi$ s'annullent sur $\prt B_{R}(0)$, $w^h_{+ }\leq \tilde L\Psi$ dans $\Gg_{R,r_{0}}$ et donc
\begin {equation}\label {est11}
\myfrac{d^2u(r,\gg_{j}(t))}{dt^2}|_{t=0}\leq \tilde L\Psi.
\end {equation}
On d\'eduit de (\ref {LB}) que $(\Gd_{S}u)_{+}(x)\to 0$ quand $\abs x\to R$. En \'ecrivant l'\'equation (\ref{equ}) en coordonn\'ees sph\'eriques, on obtient donc
\begin {equation}\label {est12}
\myfrac{\prt }{\prt r}\left(r^{N-1}\myfrac{\prt u}{\prt r}\right)\geq r^{N-1}g_{\infty}(u)+\circ (1)\quad\mbox {uniform\'ement quand }\; \abs x\to R.
\end {equation}
Clairement $u\geq z$ o\`u $z$ est la solution de 
\begin {equation}\label {est13}\left\{\BA{l}
-\Gd z+g_{\infty}(z)=0\quad\mbox {dans }\;\Gg_{R,r_{0}}\\
\,\lim_{\abs {x}\to R}z(x)=\infty\\
\;\;\;z=\min_{\prt B_{r_{0}}(0)} u \quad\mbox {sur }\;\prt B_{r_{0}}(0).
\EA\right.\end {equation}
Donc $g_{\infty}(u)\geq g_{\infty}(z)$. Comme $g_{\infty}(z)\notin L^1(\Gg_{R,r_{0}})$, on a
$$\lim_{r\to R}\myint{r_{0}}{r}g_{\infty}(u(s,\gs))s^{N-1}ds=+\infty\quad \mbox {uniform\'ement pour }\;\gs\in S^{N-1}.
$$
Ceci implique
$$\lim_{r\to R}\myfrac{\prt u}{\prt r}(r,\gs)=+\infty\quad \mbox {uniform\'ement pour }\;\gs\in S^{N-1},
$$
et donc (\ref{infde})-ii et le Lemme 1. Le th\'eor\`eme 2 en d\'ecoule.

%%%END DOCUMENT%%%%%%%%%%%%%%%

\end {document}